\documentclass[12pt]{article}
\usepackage{amssymb}

\newtheorem{theorem}{Theorem}[section]
\newtheorem{lemma}{Lemma}[section]
\newtheorem{corollary}{Corollary}[section]
\newtheorem{proposition}{Proposition}[section]
\newtheorem{definition}{Definition}[section]

\newcommand{\qed}{\hfill $\diamondsuit$\vskip .2cm}

\newcommand{\sect}[1]{\section{#1}\setcounter{equation}{0}}
\newcommand{\be}{\begin{equation}}
\newcommand{\ee}{\end{equation}}
\newcommand{\bea}{\begin{eqnarray}}
\newcommand{\eea}{\end{eqnarray}}
\newcommand{\ba}{\begin{array}}
\newcommand{\ea}{\end{array}}
\newcommand{\beas}{\begin{eqnarray*}}
\newcommand{\eeas}{\end{eqnarray*}}

\newcommand{\ZZ}{{\bf Z}}
\newcommand{\RR}{{\bf R}}
\newcommand{\NN}{{\bf N}}
\newcommand{\nn}{\nonumber}
\def\e{\varepsilon}
\def\etatil{\widetilde{\eta}}
\def\gammatil{\widetilde{\gamma}}
\def\ztil{\widetilde{z}}

\def\P{{{\cal{P}}_\rho}}
\def\v{\vskip .2cm}

\begin{document}

 \title{Transience of Second-Class Particles and 
Diffusive
     Bounds for Additive Functionals in One-Dimensional Asymmetric Exclusion 
 Processes}
 \author{Timo Sepp\"al\"ainen \  and \ Sunder Sethuraman\\
 Dept. of Mathematics\\
 Iowa State University\\
 Ames, IA \ 50011}

 \thispagestyle{empty}
 \maketitle
 \abstract{Consider a one dimensional exclusion process with
    finite-range translation-invariant jump rates with non-zero
    drift. Let the process be stationary with product Bernoulli 
invariant distribution at density $\rho$. Place a second class
particle initially at the origin. For the case  $\rho\neq
 1/2$ we show that the time spent by the second class particle  
 at the origin has finite expectation. 
 This strong transience is then used to prove that
 variances of additive functionals of local mean-zero functions are
 diffusive when $\rho\neq 1/2$. 
 As a corollary to previous work, we deduce the invariance principle 
 for these functionals. 
 The main arguments are comparisons of $H_{-1}$ norms, a large
 deviation estimate for second-class particles,  and
 a relation between occupation times of second-class
 particles and additive functional variances.
 }
 
 \v
 \thanks{$^*$ Sepp\"al\"ainen and Sethuraman are
 partially supported by NSF grants DMS-9801085 and DMS-0071504 respectively.
 \\[.15cm]
 {\sl Abbreviated title}:  Second-Class Particles and Additive
 Functionals \\[.15cm]
 {\sl AMS (1991) subject classifications}: Primary 60K35; secondary 60F05.}
 \eject

\sect{Introduction}
Informally, the simple exclusion process
updates the motion of a collection of indistinguishable random walks
on the lattice $\ZZ^d$ such that jumps to already
occupied vertices are suppressed.
These systems have had application to a wide variety of scientific
problems in physics, traffic, queuing, biology etc.  In this paper, we
exploit a connection between the diffusive behavior
of occupation times, say at a fixed
location on the lattice, and the recurrence-transience properties
of so-called ``second-class'' particles in
the exclusion model to prove results in both directions.

Briefly, we survey some of the work for second-class particles and additive 
functional
fluctuations.
The
study of the fluctuations of occupation times, or more generally
that of additive functionals, for the exclusion process was begun
by Kipnis and Varadhan \cite{KV} where they proved an invariance
principle to Brownian motion under diffusive scaling for reversible
processes in equilibrium which have finite variance.  Not all
variances
of additive
functionals are diffusive and the exact class of diffusive
additive functionals for reversible models was characterized by
Sethuraman and Xu \cite{SX}.  Notably, the occupation times in
dimensions $d=1$ and $2$ are super-diffusive, but in the appropriate
scales, $t^{3/4}$ and $\sqrt{t\log t}$ respectively,
their fluctuations were described by Kipnis \cite{Kip-flu}.  Varadhan
later generalized the Kipnis-Varadhan theorem to systems with
mean-zero jump rates \cite{V}.  Subsequently, for models
whose jump rates possess non-zero drift in $d\geq 3$, diffusive
variance bounds for all additive functionals was proved along
  with an associated invariance
principle by Sethuraman, Varadhan and Yau \cite{SVY}.  Recently, in
dimensions $d=1$ and $2$, invariance principles for some additive
functionals for models with drift
were proved provided their variances were diffusive by
Sethuraman \cite{S}.  One of the purposes of this article is to supply
the needed variance estimates to complete the story in dimension $d=1$
when $\rho\neq 1/2$
(Theorem \ref{addfctvar}).  What remains is to capture the variance
behavior for models with drift
in
dimension $d=2$ and also for $d=1$ when the density is $\rho=1/2$

Roughly, a
second-class particle in the exclusion system is a particle which
moves as a regular particle except that it also exchanges places with
regular particles which jump onto it.  In other words, the
second-class particle moves from vertex $i$ to $j$ if it jumps to an
open site at $j$ or if a particle at $j$ jumps to the position $i$.  Hence, 
the regular particle,
``first-class'' particles, do not ``see'' the second-class particle.
These particles make natural appearances in various contexts such as
in (1) the description of shocks and currents (Ferrari and Fontes 
\cite{FF}),
(2) the proofs of extremality of
some invariant measures (Saada \cite{Saada}),
and as mentioned (3) the diffusive behavior
of additive functionals.
In the third context, it is seen that the transience of a
second-class particle is equivalent to diffusive occupation-time
variance estimates.  So, in particular, in dimensions $d\geq 3$ by the
variance bounds in \cite{SVY} one concludes that the second-class
particle is transient.  One of the main results in this note is to show
that in $d=1$ for densities $\rho \neq 1/2$ when the model has drift
that the second-class
particle is also transient (Theorem \ref{secondclass}).  What is left 
open is
the recurrence-transience behavior of the particle in $d=1$ when $\rho =
1/2$ and also in $d=2$ when the system has drift.

The method of proof of the two main results, Theorems \ref{secondclass} and
\ref{addfctvar}, is to go back and forth along the bridge linking diffusive
additive functional behavior and transience of second-class particles
with the aid of two recent papers, one which gives a microscopic
variational formula for the second-class position in a specific
($K$-)exclusion model, Sepp\"al\"ainen \cite{Sepp-K},
and one which proves that diffusive variances
in one process with drift is equivalent to diffusive variances in many other
processes with drift, Sethuraman \cite{comparison}.  The strategy is to prove a
second-class particle large deviation estimate for a specific
exclusion model with drift in $d=1$ for $\rho \neq 1/2$ following from
a variational relation proved in \cite{Sepp-K}.  The
large deviation result
will imply diffusive additive functional
variance bounds for this model.  Using \cite{comparison}, we then
get
that all models with drift in $d=1$ and $\rho\neq 1/2$
have diffusive additive functional bounds.
Therefore, translating back, second-class particles in all models with
drift in $d=1$ and $\rho \neq 1/2$ are also transient.

\sect{Definitions and Results}

To state more carefully the results, we
now define more specifically
the exclusion model and the notion of a second-class particle.
Let $\Sigma = \{0,1\}^{\ZZ^d}$ be the configuration space and let
$\eta(t) \in \Sigma$ be the state of the process at time $t$.
The exclusion configuration is usefully given in terms of
occupation variables
$\eta(t) = \{\eta_i(t): i\in \ZZ^d\}$ where
$\eta_i(t)=0
{\rm \ or \ }1$
according to whether
the vertex $i\in \ZZ^d$ is empty or full at time $t$.  Let
$\{p(i,j): i,j\in \ZZ^d\}$ be the random walk or
particle transition rates.
Throughout this article
we concentrate on the translation-invariant finite-range case:
$p(i,j)=p(0,j-i)=p(j-i)$ and $p(x)=0$ for $|x|>R$ some integer
$R<\infty$.  In addition, to avoid technicalities,
we will also assume that the symmetrization
$\bar{p}(i) = (p(i) +p(-i))/2$ is irreducible.

The evolution of the system
$\eta(t)$ is Markovian.  Let $\{T_t, t\geq 0\}$ denote the
process semi-group and let
$L$ denote the infinitesimal generator.  On test functions $\phi$,
$(T_t\phi)(\eta)=E_\eta[\phi(\eta(t))]$ and
\be
\label{generator}
(L\phi)(\eta) = \sum_{i,j}
\eta_i(1-\eta_j)(\phi(\eta^{i,j})-\phi(\eta))p(j-i)
\ee
where $\eta^{i,j}$ is the ``exchanged'' configuration, $(\eta^{i,j})_i 
=\eta_j$,
$(\eta^{i,j})_j = \eta_i$ and $(\eta^{i,j})_k = \eta_k$ for $k\neq i,j$.
The transition
rate $\eta_i(1-\eta_j)p(j-i)$ for $\eta \rightarrow \eta^{i,j}$
represents the exclusion property.  The construction of the infinite
particle system follows from the Hille-Yosida theorem or by graphical
methods \cite{Lig}.

The equilibria for the exclusion system are well known.
As the exclusion model is conservative,
in that random-walk particles
are neither
destroyed nor created,
one expects a family of invariant measures indexed according
to particle density $\rho$.  In fact, let
$P_\rho$, for $\rho \in [0,1]$, be the
infinite Bernoulli product measure over $\ZZ^d$ with marginal
$P_\rho\{\eta_i =1\}=1-P_\rho\{\eta_i-0\}=\rho$.
It is shown in \cite{Lig} that $\{P_\rho:
\rho \in [0,1]\}$
are invariant for $L$.  In fact, it is proved
in \cite{Saada} that the $P_\rho$ for $\rho\in[0,1]$
are also extremal
in the convex set of invariant measures for $L$.

Let the path measure with initial distribution $P_\rho$ be given
by ${\cal{P}}_\rho$.
Let
$E_{\mu}$ be
expectation with respect to the measure $\mu$.
When the context is clear, we will denote $E_\mu$ for
$\mu =P_\rho$, or ${\cal{P}}_\rho$
as simply $E_\rho$.

We now turn to the definition of a second-class particle in the
exclusion set-up.  Consider two initial configurations $\eta$ and
$\etatil$ such that $\etatil_i = \eta_i$ for all $i\neq 0$ and
$\etatil_0=1-\eta_0=1$.  Let
$\eta(t)$ and $\etatil(t)$ be the corresponding exclusion
configurations at time $t\geq 0$.
By the
basic coupling, or ``attractive'' nature of exclusion processes (Liggett 
\cite{Lig}), we may
couple the two processes so that $\eta(t)$ and $\etatil(t)$ also differ
at exactly one vertex at any time $t\geq 0$.  Let $R(t)$ be the
position of this discrepancy, or extra particle in the $\etatil(t)$
system, at time $t$.
Then $R(t)$ describes a
second-class particle.  This can be read from the joint
generator
$\hat{L}$ of the
$(\etatil(t), R(t))$ system:
\begin{eqnarray*}
(\hat{L}\phi)(\etatil,r) &=& \sum_{i,j\neq r}
[\etatil_i(1-\etatil_j)p(j-i)](\phi(\etatil^{i,j},r)-\phi(\etatil,r))\\
&&\ \  + \sum_{k}[\etatil_{r-k}p(k) + (1-\etatil_{r-k})p(-k)]
(\phi(\etatil^{r,r-k}, r+k)-\phi(\etatil,r)).
\end{eqnarray*}
Here, the rate $\etatil_{r-k}p(k) + (1-\etatil_{r-k})p(-k)$ represents two
possible movements, namely
when a dominant particle in the $\etatil$
configuration moves to the discrepancy position $r$, and
when the
discrepancy particle jumps to an empty site at $r-k$.
The difficulty in the analysis of $R(t)$ is that it is not Markovian in general
with respect to its own history.  However, it is notable that when the
jump rate $p$ is symmetric, then the rate for $r\rightarrow r-k$
simplifies to
$\etatil_{r-k}p(k) + (1-\etatil_{r-k})p(-k)= p(k)$ so that in this case
$R(t)$ is a bona fide symmetric random-walk.

We now describe the connection between occupation-times and
second-class particles known in the folklore.  Consider the exclusion system in equilibrium
$P_\rho$.  Let $f(\eta)= \eta_0 - \rho$ be the centered occupation
function.  Let us compute the variance of the occupation time of the
origin up to time $t$:
$$E_\rho[\big ( \int_0^t (\eta_0(s)-\rho) ds\big )^2]
= 2\int_0^t (t-s)E_\rho[f(\eta(s))f(\eta(0))]ds.$$
We may expand the kernel further,
\begin{eqnarray*}
E_\rho[f(\eta(s))f(\eta(0))] &=&
E_\rho[\eta_0(s)\eta_0(0)]-\rho^2\\
&=&\rho\{\P[\eta_0(s)=1|\eta_0(0)=1]-\P[\eta_0(s)=1]\}\\
&=&\rho(1-\rho)\{\P[\eta_0(s)=1|\eta_0(0)=1] -
\P[\eta_0(s)=1|\eta_0(0)=0]\}.
\end{eqnarray*}
To rewrite the
last difference further, we couple the initial measures
$P_\rho(\cdot|\eta_0=1)$ and $P_\rho(\cdot|\eta_0=0)$ through the
basic coupling so that the two
systems differ in only one position at any later time.
This discrepancy
position is of course the second-class position in a sea of
regular particles distributed initially according to
$P_\rho(\cdot|\eta_0=0)$.  The last line, therefore,
under the coupling measure $\bar{\P}$ is
restated as $\rho(1-\rho)\bar{\P}[R(s)=0]$.
Evidently then the occupation-time variance satisfies the
following relation with the expected occupation time of the second-class 
particle:
\be
\label{relation}
\lim_{t\rightarrow \infty}
\frac{1}{t}E_\rho[\big ( \int_0^t (\eta_0(s)-\rho) ds\big )^2]\ = \
2\rho(1-\rho)\int_0^\infty \bar{\P}[R(s)=0]ds.
\ee

This calculation motivates the following definition.
\begin{definition}
The second-class particle is $\P$-{\it recurrent} or {\it
   transient at } $0$ if respectively
$$\int_0^\infty \bar{\P}[R(s)=0]ds = \infty \ \ {\rm or \ \ }
<\infty.$$
\end{definition}

One of the main results in this paper is the following.
\begin{theorem}
\label{secondclass}
For simple exclusion processes in $d=1$ with finite-range
translation-invariant jump rates $p$ which have non-zero drift,
$\sum_i ip(i)\neq 0$, the second-class particle is $\P$-transient
at $0$ when the equilibrium density $\rho \neq 1/2$.
\end{theorem}

Central to the proof of this theorem is the
 large deviation estimate below, of interest in its
own right.
\begin{theorem}
\label{secondclassestimate}
For the simple exclusion process in $d=1$ with totally asymmetric
nearest-neighbor translation-invariant
jump rates $p$ such that $p(1)=1$ and $p(i)=0$ for all $i\neq 1$,
 there exist constants $A=A(\e,\rho)$ and 
$C=C(\e,\rho)>0$ such that for all $t>0$, 
$$\bar{\P}[ \ |R(t)-(1-2\rho)t|>\e t \ ] \ \leq \ A
e^{-Ct}.$$
\end{theorem}

The coefficient $v(\rho,p)=(1-2\rho)\sum_i ip(i)$ (which reduces to 
$1-2\rho$ in
the above theorem) is the limiting velocity of the second-class
particle. It was proved by Ferrari \cite{Ferrari} that
\be
\label{transiencelimit}
\lim_{t\rightarrow\infty} \frac{R(t)}{t} = v(\rho,p) \ \ {\rm a.s. \
   } (\bar{\P})
\ee
for nearest-neighbor translation-invariant processes.
Our proof of Theorem \ref{secondclassestimate} is independent
of Ferrari's law of large numbers 
 and so  gives
another proof of (\ref{transiencelimit}).
At $\rho\ne 1/2$, the limit (\ref{transiencelimit})
 gives transience in the usual sense,
\be
\label{usualtransience}
\bar{\P}[ \ R(\cdot) \ {\rm visits \ }0 \ {\rm finitely \ many \ times
   } \ ] =1.
\ee
   But, unfortunately, we could not directly
convert (\ref{usualtransience}) to
   $\P$-transience at the origin.
The converse holds however.  As the jump rates of
   $R(t)$ are bounded,
$$\sum_k[\eta_{r-k}p(k) + (1-\eta_{r-k})p(-k)] \ \leq \
   \sum_k[p(k)+p(-k)] \ <\infty$$
uniformly in the environment, we can couple an independent exponential
r.v. $U$ having intensity
$\sum_k[p(k)+p(-k)]$ with the jump time variable
$\tau$ each time the second-class particle visits the origin.
Therefore,
$$\int_0^\infty\bar{\P}[R(t)=0]dt \ \geq \
\frac{1}{\sum_k[p(k)+p(-k)]} \bar{E_\rho}[\# {\rm \ visits} \ R(\cdot) {\rm
   \ makes \ to \ the \ origin}]$$
and so, $\P$-transience implies usual transience.
A similar argument shows that recurrence, in the usual sense,
$$
\bar{\P}[\ \exists \{t_n\}, t_n\uparrow\infty {\rm \ s.t. \ }
R(t_n)=0\ ]=1,$$
implies $\P$-recurrence.

Other related results on the second-class particle are that $R(t)$ is
$\P$-recurrent in $d=1,2$ when the jump rate is mean-zero, $\sum_i
ip(i)=0$.  Also, in $d\geq 3$, $R(t)$ is $\P$-transient no matter
what the jump rate $p$ is.  See section 6 \cite{S} for details and extensions.

The open cases are when
the jump rate has non-zero drift in $d=2$, and in $d=1$ with $\rho=1/2$.
The latter situation is quite tantalizing
as there seems to be intuition for both $\P$-recurrence and
$\P$-transience.  On the one hand, it should be $\P$-recurrent due to
the analogy with random walk with zero velocity.
But, on the other hand, there could be
a remnant of the behavior of second class particles in a rarefaction
fan which leads to transient behavior.
That is, it might be possible that $R(t)$ flips a fair
coin and on the basis of the toss would end up eventually
exclusively on the left or right of the
origin.  See Ferrari and Kipnis \cite{FK} for details about the rarefaction 
fan behavior.

\vskip .5cm
We now turn to our results for the diffusive behavior of additive
functionals.  Let $f$ be a local mean-zero function, that is, a function
which depends
only on a finite number of coordinates, and $E_\rho[f]=0$.
Define
$$A_f(t) = \int_0^t f(\eta(s))ds$$
as the additive functional of $f$ up to time $t$.  Denote the variance
of $A_f(t)$ as $\sigma^2_t(f,\rho,p)= E_\rho[(A_f(t))^2]$ and denote
also the limiting variance, if it exists,
$$\sigma^2(f,\rho,p) = \lim_{t\rightarrow \infty} t^{-1}
\sigma^2_t(f,\rho,p).$$

\begin{theorem}
\label{addfctvar}
For simple exclusion processes in $d=1$ with finite-range
translation-invariant jump rates $p$ with non-zero drift, $\sum_i
ip(i) \neq 0$, and density $\rho\neq 1/2$, we have, for any local mean-zero 
function $f$, that
$\sigma^2(f,\rho,p)<\infty$.
\end{theorem}

The finiteness of the limiting variances give the following corollary.
\begin{corollary}
\label{invarprin}
In the case of Theorem \ref{addfctvar},
we have the weak convergence to
Brownian motion $B$ in $C[0,\infty)$,
\be
\label{invariance}
\lim_{\alpha\rightarrow\infty}
\frac{1}{\sqrt{\alpha}} \int_0^{\alpha t} f(\eta(s))ds
= B(\sigma^2(f,\rho,p)t).
\ee
\end{corollary}

The limiting variance $\sigma^2(f,\rho,p)$, or diffusion coefficient, has 
been itself an
object of much attention.  Only in a few specific cases has it been
explicitly computed, and for the most part almost all the 
theoretical
work has concentrated on existence proofs.  When $p$ is
symmetric in all $d\geq 1$, the coefficient is known to exist and is positive
for all non-constant local functions $f$, and in
fact half the diffusion coefficient is also the square of the $H_{-1}$ norm
of $f$,
$\sigma^2(f,\rho,p)/2 = \|f\|^2_{-1}(\rho,p)$ \cite{KV}.
It was proved later, for $p$ symmetric,
that $\|f\|_{-1}(\rho,p)<\infty$ if and only if
$d^n/d\theta^n E_\theta[f]|_{\theta = \rho} = 0$ for $n=0,1,2$ in
$d=1$, $n=0,1$ in $d=2$, and $n=0$ in $d\geq 3$ \cite{SX}.
In particular, as noted in the introduction,
the occupation function $f(\eta)=\eta_0-\rho$ is not in
$H_{-1}$ for symmetric $p$ in dimensions $d=1$ and $2$, the correct
orders in $d=1,2$ being $\sigma^2_t(\eta_0-\rho,\rho,p) \sim t^{3/2}$
and $t\log t$ \cite{Kip-flu}.  For asymmetric but mean-zero
$p$, it is shown, for non-constant local $f$, that 
$0<\sigma^2(f,\rho,p)<\infty$ if and only if $f\in
H_{-1}(\rho,\bar{p})$
\cite{V}, \cite{S}.  For asymmetric $p$ with non-zero drift, it is proved that
$0<\sigma^2(f,\rho,p)<\infty$ exists in $d\geq 3$ for all
(non-constant) mean-zero
local functions $f$ \cite{SVY},
\cite{S}--the same condition as for symmetric $p$ in $d\geq 3$.
In $d=1,2$, when $p$ has non-zero drift, it was shown that
$0<\sigma^2(f,\rho,p)$ exists when $f$ is non-constant and
increasing; existence and finiteness was also shown
when
$f=f_+ - f_-$, the difference of two increasing mean-zero functions whose
variances are finite, $\sigma^2(f_\pm,\rho,p)<\infty$ \cite{S}.
Notably, it was not shown in \cite{S} when
$0<\sigma^2(f,\rho,p)<\infty$ occurs in general.
However, it was
shown in \cite{comparison}
that for every $p$ there is a nearest-neighbor jump rate $p'$ such
that, if both $\sigma^2(f,\rho,p)$ and $\sigma^2(f,\rho,p')$
exist, then
$\sigma^2(f,\rho,p)$ is bounded  (positive) if and only if
$\sigma^2(f,\rho, p')$ is bounded  (positive).
The contribution of Theorem \ref{addfctvar}
is the statement 
 $\sigma^2(f,\rho,p)<\infty$ for local mean-zero $f$
in the case 
$d=1$ and $\rho\neq 1/2$.  The open cases left are in $d=2$ and in
$d=1$ for $\rho = 1/2$, with non-zero drift.

Existence and finiteness
of the diffusion coefficient  is half the
question, the other half being ``does a central limit theorem hold?''  The
answer is basically ``yes.''  Kipnis and Varadhan established the
invariance principle (\ref{invariance}) for $p$ symmetric when $f\in H_{-1}$
 by martingale approximation
\cite{KV}.  When $p$ is
asymmetric but mean-zero, Varadhan generalized this method
and proved the invariance principle for $f$ such that
$\sigma^2(f,\rho,\bar{p})<\infty$
\cite{V}.  When $p$ is asymmetric with non-zero drift, the invariance
principle was established for all mean-zero local
functions $f$ in $d\geq 3$
\cite{SVY}.  In $d=1,2$ when $p$ has non-zero drift, the invariance
principle was proved for increasing meanzero $f\in L^2(P_\rho)$ such
that $\sigma^2(f,\rho,p)<\infty$ and also for $f=f_+-f_-$, the
difference of local increasing mean-zero 
functions whose variances are finite,
$\sigma^2(f_\pm,\rho,p)<\infty$, through techniques with associated
r.v.'s \cite{S}.
In Corollary \ref{invarprin}, we extend this result to all mean-zero
local $f$ in $d=1$
when $\rho\neq 1/2$.  Again, the open cases revolve around $d=2$ and
$d=1$ when $\rho = 1/2$, when the jump rate $p$ has non-zero drift.




\vskip .5cm

The plan of the paper is the following.  In section 3, we prove the
hard estimate Theorem \ref{secondclassestimate}.  In section 4, we prove 
the rest of the
results,
Theorems \ref{secondclass} and \ref{addfctvar}, and Corollary \ref{invarprin}
in short succession.

\section{Proof of large deviation bound Theorem \ref{secondclassestimate}}

Now we restrict ourselves to the totally asymmetric, nearest neighbor
case, so $p(1)=1$ and $p(i)=0$ for $i\ne 1$.
If $\rho=0$ or $1$, $R(t)$ is a Poisson process
and the desired
estimate is trivial. So we assume
$\rho\in(0,1)$.
To prove Theorem \ref{secondclassestimate} we turn to the variational coupling
representation of the totally asymmetric exclusion process.

\subsection{The second class particle in the variational
coupling}

To describe this coupling, we will need to recall some details of the
  graphical construction of the exclusion model.
We perform this construction
of the exclusion process $\eta(t)$
in terms of
   a collection
  $\{D_i:i\in\ZZ\}$ of mutually
independent rate 1 Poisson jump time processes on the time line
$(0,\infty)$.  Let  $(\Omega,{\cal F}, P)$ denote a
  probability
space on which the $\{D_i\}$ are defined, and
independently of them the initial configuration
$(\eta_i(0):i\in\ZZ)$ of the exclusion process.
In the construction we represent $\eta(t)$
  in terms of
``current particles.'' These form
a process $z(t)= (z_i(t):i\in{\ZZ})$
  of labeled particles that move on  ${\ZZ}$
subject to the constraint
\be
0\le z_{i+1}(t)-z_i(t)\le 1\quad\mbox{for all $i\in\ZZ$ and $t\ge 0$.}
\label{tz1}
\ee
   In the graphical construction,
  $z_i$ attempts to jump one step to the {\it left} at
epochs of $D_i$. If the execution of the jump would
produce a configuration that
violates (\ref{tz1}), the jump is suppressed.
  We can summarize the jump rule
like this:
\bea
&&\mbox{If $t$ is an epoch of $D_i$, then}\nn\\
&&z_i(t)=\max\{ z_i(t-)-1, z_{i-1}(t-), z_{i+1}(t-)-1\}.
\label{tzjumprule}
\eea

We arrange things so that
  $\eta$
gives the increments of $z$.
Given the initial configuration
  $\{\eta_i(0)\}$, the initial configuration
$\{z_i(0)\}$ is defined on $\Omega$ by
\be
z_0(0)= 0\,,\
z_i(0)=\sum_{1\le j\le i}\eta_j(0)\ \mbox{ for $i>0$, and }
\
z_i(0)=-\sum_{i< j\le 0}\eta_j(0)\ \mbox{ for $i<0$.}
\label{tinitialzeta}
\ee
The choice $z_0(0)=0$ is merely a convenient normalization.
Any random choice independent of $\{\eta_i(0)\}$ and
  $\{D_i\}$ would do.

We can construct the process $z(t)$ by applying
the jump  rule (\ref{tzjumprule}) inductively to jump times,
  once we exclude
  an exceptional null set of ``bad'' realizations of
$\{D_i\}$.
We always
assume that the realization
$\{D_i\}$ satisfies these requirements:
\bea
&&\mbox{(i) There are no simultaneous jump attempts.}
\nn\\
&&\mbox{(ii) Each
$D_i$ has only finitely may epochs in every bounded time
interval.} \nn\\
&&\mbox{(iii) Given any $t_1>0$, there are arbitrarily faraway
indices $i_0<<0<<i_1$}\nn\\
&&\mbox{such that $D_{i_0}$
and $D_{i_1}$  have no epochs
in the time interval $[0,t_1]$.}
\nn
\eea
These properties are satisfied
almost surely, so the evolution $z(t)$, $0\le t<\infty$,
  is well-defined for
  almost every realization of $\{D_i\}$. Then the
process $\eta(t)$ is defined for $t>0$ by
\be
\eta_i(t)=z_i(t)-z_{i-1}(t)\,.
\label{tetaz}
\ee
It should be clear that
$\eta(t)$ operates as an exclusion process with jump
probabilities $p(1)=1$ and $p(i)=0$ for $i\ne 0$. 
The $z$-process
represents the current of $\eta$, for
$z_i(0)-z_i(t)$ equals the number of exclusion particles that have
jumped across the bond $(i,i+1)$ during the time interval
$(0,t]$.

For the variational coupling we construct
  a family $\{w^k:k\in{\ZZ}\}$ of auxiliary processes
on the space  $\Omega$. Each
$w^k(t)=(w^k_i(t):i\in{\ZZ})$ is a process
of the same type as $z(t)$. 
The initial configuration $w^k(0)$ depends on the initial
position $z_k(0)$:
\bea
w^k_i(0)=\left\{ \begin{array}{ll}
z_k(0)\,, &i\ge 0\\
z_k(0)+i\,, &i<0.
\end{array}
\right.
\label{twinit}
\eea
The processes $\{w^k\}$ are coupled to each other
and to $z$ through the
Poisson processes $\{D_i\}$. However,
the jump rule for $w^k$ includes a translation
of the index:
\be
\mbox{at epochs $t$ of $D_{i+k}$, }\
w^k_i(t)=\max\{ w^k_i(t-)-1, w^k_{i-1}(t-), w^k_{i+1}(t-)-1\}.
\label{twjumprule}
\ee
The 
increments process $w^k_i(t)-w^k_{i-1}(t)$
 represents an exclusion
process where initially the lattice is full from
site $k$ to the left, and empty from site $k+1$
to the right.
The point of introducing the processes
$\{w^k\}$ lies in this ``variational coupling''
lemma:

\begin{lemma} For all $i\in{\ZZ}$ and $t\ge 0$,
\be
z_i(t)=\sup_{k\in{\ZZ}} w^k_{i-k}(t)\quad\mbox{a.s.}
\label{tzsupw}
\ee
\label{tvarcouplm}
\end{lemma}

This lemma is proved by induction on jump times, assuming
properties (i)--(iii) above for $\{D_i\}$. For details,
  see Lemma 4.2 in Sepp\"al\"ainen \cite{Se2}.

For Theorem \ref{secondclassestimate}
we need deviation bounds
for the processes $w^k$. For this we
decompose $w^k$ into a sum of the initial position
defined by (\ref{twinit}) and the increment determined
by  the Poisson  processes through (\ref{twjumprule}).
To this end, define a family of processes $\{\xi^k\}$ by
$$\mbox{$\xi^k_i(t)=z_k(0)-w^k_i(t)$ for $i\in\ZZ$, $t\ge 0$.}$$
The process
$\xi^k$ does not depend on $z_k(0)$, and depends on
the superscript $k$ only through a translation of the
$i$-index of the Poisson processes
$\{D_i\}$.
  Initially
\bea
\xi^k_i(0)=\left\{ \begin{array}{ll}
0\,, &i\ge 0\\
-i\,, &i<0.
\end{array}
\right.
\label{txiinit}
\eea

Dynamically, at epochs $t$ of $D_{i+k}$,
$$
\xi^k_i(t)=\min\{ \xi^k_i(t-)+1, \xi^k_{i-1}(t-), \xi^k_{i+1}(t-)+1\}.
$$
We think of $\xi^k$ as a growth model on the upper half plane,
so that $\xi^k_i$ gives the height of the interface above
site $i$.
It can be equivalently defined by specifying that each $\xi^k_i$
advances independently at rate 1, provided   these inequalities
are preserved:
\be
\xi^k_i(t)\le \xi^k_{i-1}(t)\qquad\mbox{and}\qquad
\xi^k_i(t)\le \xi^k_{i+1}(t)+1.
\label{txiineq}
\ee

In terms of $\xi$, (\ref{tzsupw}) can be expressed as
\be
z_i(t)=\sup_{k\in{\ZZ}} \{z_k(0)-\xi^k_{i-k}(t)\}.
\label{tzsupxi}
\ee

Next we include the
  second class particle $R(t)$  in the variational
coupling picture. Recall the definition of $R(t)$
as the location of the unique discrepancy between two
processes  $\eta$ and $\etatil$ that initially agree
everywhere except at $R(0)$, where
$\etatil_{R(0)}(0)=1-\eta_{R(0)}(0)=1$.  [Earlier
we took $R(0)=0$ but that is not necessary for what follows here.]

We define $\eta$ and $\etatil$ by (\ref{tetaz}),  in terms
of processes $z$ and $\ztil$ that initially satisfy
$$
\mbox{$\ztil_i(0)=z_i(0)$ for $i\le R(0)-1$
and $\ztil_i(0)=z_i(0)+1$ for $i\ge R(0)$. }
$$
It may happen that $\ztil_0(0)\ne 0$, but that is of
no consequence. We make the
processes $z$ and $\ztil$ obey the same Poisson
processes $\{D_i\}$ through the jump rule
(\ref{tzjumprule}), so this is the basic coupling.
One can prove that at all times $t\ge 0$
  there is a unique discrepancy marked
by  $R(t)$:
\be
\mbox{$\ztil_i(t)=z_i(t)$ for $i\le R(t)-1$
and $\ztil_i(t)=z_i(t)+1$ for $i\ge R(t)$. }
\label{tXdef2}
\ee
Using (\ref{tXdef2}), we prove a variational
representation for $R(t)$.

\begin{proposition} Almost surely, for all $t\ge 0$,
\be
R(t)=\inf\{ i\in{\ZZ}: \mbox{$z_i(t)=z_k(0)-\xi^k_{i-k}(t) $ for
some $k\ge R(0)$}\}.
\label{tXinf}
\ee
\label{tXthm}
\end{proposition}

{\it Proof.}  The claim (\ref{tXinf}) will follow from proving
\be
\mbox{if $i<R(t)$, then $z_i(t)>z_k(0)-\xi^k_{i-k}(t) $ for
all $k\ge R(0)$,}
\label{tXinf1}
\ee
and
\be
\mbox{if $i\ge R(t)$, then $z_i(t)=z_k(0)-\xi^k_{i-k}(t) $ for
some $k\ge R(0)$.}
\label{tXinf2}
\ee

To contradict (\ref{tXinf1}), suppose  $i<R(t)$ and
$z_i(t)=z_k(0)-\xi^k_{i-k}(t)
$
for
some $k\ge R(0)$. Then by (\ref{tXdef2})
$$\ztil_i(t)=\ztil_k(0)-\xi^k_{i-k}(t)-1$$
which contradicts the variational formula (\ref{tzsupxi})
for process $\ztil$. This contradiction proves (\ref{tXinf1}).

To prove  (\ref{tXinf2}), let $i\ge R(t)$. Suppose
that for some  $k< R(0)$,
$\ztil_i(t)=\ztil_k(0)-\xi^k_{i-k}(t)$. Then by  (\ref{tXdef2}),
$$z_i(t)+1=z_k(0)-\xi^k_{i-k}(t).
$$
This contradicts (\ref{tzsupxi}), so it must be that
$\ztil_i(t)=\ztil_k(0)-\xi^k_{i-k}(t)$ for some  $k\ge R(0)$.
Again by  (\ref{tXdef2}), this implies
$z_i(t)=z_k(0)-\xi^k_{i-k}(t)
$. This proves (\ref{tXinf2}).
\qed

\subsection{Auxiliary results}

As mentioned, all the processes $\xi^k$ have the same
distribution, because the effect of the superscript
$k$ is only to translate the index of the Poisson jump
time processes $\{D_i\}$. Let us write $\xi$ to
simultaneously denote any one of them. A law of
large numbers is given by
\be
\lim_{t\to\infty}t^{-1} \xi_{[tx]}(t)=g(x)
\quad\mbox{almost surely,}
\label{tllnxi}
\ee
where $g$ is defined by
\be
g(x)=\left\{\begin{array}{ll}
-x, &x<-1\\
(1/4)(1-x)^2, &-1\le x\le 1\\
0, &x\ge 1.
\end{array}
\right.
\label{tgdef}
\ee
This result goes back to Rost \cite{Ro}.
For $\xi$ we have these large deviation bounds.

\begin{proposition} Let $x\in\RR$ and $\e>0$. Then there exists
a finite positive
constant $C$ such that for all $t>0$,
\be
P\left( \xi_{[tx]}(t)\ge tg(x)+t\e\right)\le \exp(-Ct^2)
\label{tldp2ut}
\ee
and
\be
P\left( \xi_{[tx]}(t)\le tg(x)-t\e\right)\le \exp(-Ct).
\label{tldp2lt}
\ee
\label{tldp2}
\end{proposition}

{\it Proof.} We can infer this proposition from the results
in Sepp\"al\"ainen \cite{Se1} via a simple mapping of the lattice.
The first step is to convert $\xi$ into a last-passage model.
Define the passage times $L_{i,j}$ by
\be
L_{i,j}=\inf\{ t\ge 0: \xi_i(t)\ge j\}
\label{tLinf}
\ee
for $i\in\ZZ$ and $j\ge \max\{0,-i\}$. From the rules of
$\xi$ we infer the boundary conditions
$
\mbox{$L_{-i,i}=L_{i,0}=0$ for $i\ge 0$,}
$
and the equation
$$
L_{i,j}=\max\left\{L_{i-1,j}\,,\,L_{i,j-1}\,,\,L_{i+1,j-1}\right\}
+ Y_{i,j} \ \mbox{ for $j>\max\{0,-i\}$,}
$$
where $Y_{i,j}$ is a rate 1 exponential waiting time, independent
of the $L$-variables in braces on the right-hand side. Applying
this relation inductively leads to
\be
L_{i,j}=\max_\pi\sum_{{\bf u}\in\pi}Y_{\bf u},
\label{tLmax}
\ee
where the maximum is over lattice paths $\pi=\{ (0,1)=(i_0,j_0),(i_1,j_1),
\ldots,(i_n,j_n)=(i,j)\}$ that take
three types of steps:
$$
\mbox{$(i_{m+1},j_{m+1})-(i_m,j_m)=(-1,1)$, $(0,1)$, or $(1,0)$
for each $m$.}
$$
Eqs.\  (\ref{tLinf}) and (\ref{tLmax}) give two
different constructions
of the process  $\{L_{i,j}\}$:  in (\ref{tLinf})
  in terms of the Poisson processes $\{D_i\}$, but in
  (\ref{tLmax})  in terms of the i.i.d.\ exponential random
variables $\{Y_{i,j}\}$. The last-passage formulation
  (\ref{tLmax}) is convenient for large deviation analysis.
Corresponding to (\ref{tllnxi})--(\ref{tgdef}) we have
the strong law of large numbers
\be
\lim_{t\to\infty} t^{-1}L_{[tx],[ty]}=\gamma(x,y)\equiv
\left( \sqrt{x+y\,}+\sqrt{y}\,\right)^2
\ \mbox{for $y> 0\vee (-x)$.}
\label{tllnL}
\ee
The connection between the limits in (\ref{tllnxi}) and
(\ref{tllnL}) is, naturally enough, that the limiting
interface $g$ is a level curve of the limiting passage time:
$\gamma(x, g(x))=1$ for $-1\le x\le 1$.

Now (\ref{tldp2ut})--(\ref{tldp2lt}) will follow from proving
\be
P\left( L_{[tx],[ty]}\le t\gamma(x,y)-t\e\right)\le \exp(-Ct^2)
\label{tldp3lt}
\ee
and
\be
P\left(L_{[tx],[ty]}\ge t\gamma(x,y) +t\e\right)\le \exp(-Ct).
\label{tldp3ut}
\ee
This is exactly what is proved in \cite{Se1} for a passage-time
process $\{T_{k,l}:(k,l)\in\NN^2\}$
  that is essentially the same as $\{L_{i,j}: i\in\ZZ, j\ge 1+
(0\vee(-i))\}$. Here $\NN=\{1,2,3,\ldots\}$ is the set of
natural numbers. 
To define $T_{k,l}$, let $\{W_{k,l}:(k,l)\in\NN^2\}$ be i.i.d.\ exponential
mean 1 random variables, and set
\be
T_{k,l}=\max_\sigma\sum_{{\bf u}\in\sigma}W_{\bf u},
\label{tTmax}
\ee
where the maximum is over lattice paths $\sigma
=\{ (1,1)=(k_0,l_0),(k_1,l_1),
\ldots,(k_n,l_n)=(k,l)\}$ in $\NN^2$ that take
only up-right steps:
$$
\mbox{$(k_{m+1},l_{m+1})-(k_m,l_m)= (0,1)$ or $(1,0)$
for each $m$.}
$$

To find the correspondence between
$L_{i,j}$ and $T_{k,l}$,  observe first that the optimal
path $\pi$ in (\ref{tLmax}) never uses a $(0,1)$-step
because such a step can be replaced by a $(-1,1)$-step
followed by a $(1,0)$-step. So in (\ref{tLmax}),
let us consider only paths $\pi$ with steps
  $(-1,1)$ and $(1,0)$. Let $\Psi$ be the bijective
map from $\{(i,j)\in\ZZ\times\NN: j\ge 1+(0\vee(-i))\}$
onto $\NN^2$ given by $\Psi(i,j)=(i+j,j)$. Then
under $\Psi$ and $\Psi^{-1}$ the paths $\pi$ and $\sigma$ map onto
each other. If in (\ref{tTmax})
  we take $W_{k,l}=Y_{\Psi^{-1}(k,l)}$,
then $L_{i,j}=T_{\Psi(i,j)}$. Combining
Theorems 7, 8 and 10 in \cite{Se1} gives the estimates
(\ref{tldp3lt})--(\ref{tldp3ut}) with $L$ replaced
by $T$, and with $\gamma(x,y)$
replaced by the limit $\gammatil(x,y)=
\left( \sqrt{x}+\sqrt{y}\,\right)^2$ of $t^{-1}T_{[tx],[ty]}$.
Via $\Psi$ the estimates for $T$ become exactly
(\ref{tldp3lt})-(\ref{tldp3ut}) for $L$.
\qed

{\it Remark.} From the work of Johansson \cite{Jo} one can
get better estimates for the distribution of $T_{[tx],[ty]}$,
but this is not needed for our purposes.

\begin{lemma} Fix $k<l$. Then almost surely
$
  \xi^k_{i-k}(t) \le \xi^l_{i-l}(t)
$
for all $i\in\ZZ$ and $t\ge 0$.
\label{torderlm}
\end{lemma}

{\it Proof.} The statement is valid at time $t=0$ by
(\ref{txiinit}), and consequently valid
at all $t\ge 0$ because the coupling preserves
  ordering. (Note that both $ \xi^k_{i-k}$ and $\xi^l_{i-l}$
jump at epochs of $D_i$.)
\qed

\begin{lemma} Let   $a<x-1<x+1<b$. Then there exists a finite
constant $C\in(0,\infty)$ such that for all $t>0$,
$$
P\left( z_{[tx]}(t)\ne \max_{[ta]\le k\le [tb]}
\left\{ z_k(0)-\xi^k_{[tx]-k}(t)\right\} \right) \le \exp(-Ct).
$$
\label{tindexlm}
\end{lemma}

{\it Proof.} The initial arrangement (\ref{txiinit}) and the
constraint (\ref{txiineq}) together imply that for
$j>0$, the first jump of $\xi_j$ cannot happen
before the first jump of $\xi_{j-1}$, and correspondingly
for $j<0$. Since waiting times are exponential,
it follows that
for any $j\in\ZZ$,
the time when $\xi_j$ first jumps is distributed
as the sum of $|j|+1$
i.i.d.\ rate 1 exponential random variables.
Since $j=[tx]-[ta]$ and $j=[tx]-[tb]$ both satisfy
$|j|\ge t(1+\delta)$ for some $\delta>0$ for large enough $t$,
  standard i.i.d.\ large deviation bounds give
$$
P\left(
\xi^{[ta]}_{[tx]-[ta]}(t)=0\ \mbox{ and }\
\xi^{[tb]}_{[tx]-[tb]}(t)=[tb]-[tx]
\right)\ge 1-e^{-Ct}.
$$

To prove the lemma, it remains to check that
on the event
$$
\left\{
\xi^{[ta]}_{[tx]-[ta]}(t)=0\ \mbox{ and }\
\xi^{[tb]}_{[tx]-[tb]}(t)=[tb]-[tx]
\right\}
$$
we have
$$
z_{[tx]}(t)= \max_{[ta]\le k\le [tb]}
\left\{ z_k(0)-\xi^k_{[tx]-k}(t)\right\}.
$$
This follows from the constraints on $\xi$ and
from (\ref{tz1}): for $k<[ta]$,
$$
z_k(0)-\xi^k_{[tx]-k}(t)= z_k(0)\le z_{[ta]}(0)
=z_{[ta]}(0)-\xi^{[ta]}_{[tx]-[ta]}(t),
$$
and  for $k>[tb]$,
$$
z_k(0)-\xi^k_{[tx]-k}(t)= z_k(0)-(k-[tx])\le z_{[tb]}(0)
-([tb]-[tx])
=z_{[ta]}(0)-\xi^{[tb]}_{[tx]-[tb]}(t).
$$
This shows that indices outside the range $[ta]\le k\le [tb]$
cannot alter the supremum.
\qed

\subsection{Proof of
Theorem \ref{secondclassestimate}}

Now return to the setting of Theorem \ref{secondclassestimate}.
Fix $\rho\in(0,1)$.  Place initially a second
class particle at the origin, so $R(0)= 0$ and
$\eta_0(0)=0$ with probability 1. For
  $i\ne 0$ the initial occupation
variables $\eta_i(0)$ are i.i.d.\ with $P(\eta_i(0)=1)=\rho$.
And then the initial configuration $\{z_i(0)\}$ is defined
by (\ref{tinitialzeta}), with $z_{-1}(0)=z_0(0)=0$.
Theorem \ref{secondclassestimate} is proved
  in two steps: the lower tail and upper tail
estimate. We let $A$ and $C$ denote finite positive constants
whose values may change from one inequality to the
next but never depend on $t$.

\subsubsection{Lower tail bound}
Let $r=1-2\rho$ throughout the proof, and $\e>0$. In this
subsection we prove
\be
P\left( R(t)\le tr-t\e\right) \le A\exp(-Ct).
\label{tldplt}
\ee
By statement (\ref{tXinf2}) and Lemma \ref{tindexlm} applied
to $x=r-\e$, we get
\beas
&&P\left( R(t)\le tr-t\e\right) \\
&\le&
P\left( \mbox{$z_{[t(r-\e)]}(t)= z_k(0)-\xi^k_{[t(r-\e)]-k}(t)$
for some $k\ge 0$}
\right) \\
&\le& e^{-Ct} +
P\left( \mbox{$z_{[t(r-\e)]}(t)= z_k(0)-\xi^k_{[t(r-\e)]-k}(t)$
for some $0\le k\le bt$}
\right).
\eeas
Check that $\rho y-g(r-\e-y)$ is strictly decreasing
for $y\ge-\e$. Choose $\delta>0$ so that
\be
\mbox{$-\rho\e-g(r) \ge \rho y-g(r-\e-y) +5\delta$ for all $y\ge 0$.}
\label{tdeltadef1}
\ee
Choose a partition $0=y_0<y_1<\cdots<y_n=b$ of $[0,b]$ so that
$y_{i+1}-y_i<\delta$ for all $i$. Then it follows that
\be
\mbox{$-t\rho\e-tg(r) \ge t\rho y_i-tg(r-\e-y_{i-1}) +4\delta t$
for all $t> 0$, $1\le i\le n$.}
\label{tineq1}
\ee
Note that, by the ordering of the $z_k$'s
and by  Lemma \ref{torderlm},
\be
\mbox{$z_k(0)-\xi^k_{[t(r-\e)]-k}(t)\le
  z_{[ty_i]}(0)-\xi^{[ty_{i-1}]}_{[t(r-\e)]-[ty_{i-1}]}(t)$
for all $[ty_{i-1}]\le k\le[ty_i]$. }
\label{tineq2}
\ee
Continue the estimation from above. First use (\ref{tineq2}),
and note that by (\ref{tzsupxi}),
$z_{[t(r-\e)]}(t)\ge z_j(0)-\xi^j_{[t(r-\e)]-j}(t)$
for $j=-[t\e]$. Then use (\ref{tineq1}). 
\beas
&&P\left( R(t)\le tr-t\e\right) \\
&\le& e^{-Ct} +
\sum_{i=1}^n
P\left( \mbox{$z_{[t(r-\e)]}(t)= z_k(0)-\xi^k_{[t(r-\e)]-k}(t)$
for some $ty_{i-1}\le k\le ty_i$}
\right)\\
&\le& e^{-Ct} +
\sum_{i=1}^n
P\left( z_{[t(r-\e)]}(t)
\le z_{[ty_i]}(0)-\xi^{[ty_{i-1}]}_{[t(r-\e)]-[ty_{i-1}]}(t)
\right)\\
&\le& e^{-Ct} +
\sum_{i=1}^n
P\left( z_{-[t\e]}(0)-\xi^{-[t\e]}_{[t(r-\e)]+[t\e]}(t)
\le z_{[ty_i]}(0)-\xi^{[ty_{i-1}]}_{[t(r-\e)]-[ty_{i-1}]}(t)
\right)\\
&\le& e^{-Ct} +
\sum_{i=1}^n \left\{
P\left( z_{-[t\e]}(0)\le -t\rho\e-\delta t\right)
+ P\left(
\xi^{-[t\e]}_{[t(r-\e)]+[t\e]}(t)\ge tg(r)+\delta t\right)\right.\\
&&\qquad\left.
  + P\left(  z_{[ty_i]}(0)\ge  t\rho y_i+\delta t\right)
+ P\left(\xi^{[ty_{i-1}]}_{[t(r-\e)]-[ty_{i-1}]}(t) \le
tg(r-\e-y_{i-1}) -\delta t\right)
\right\}\\
&\le& Ae^{-Ct}.
\eeas
In the last step we use Proposition \ref{tldp2} for the probabilities
involving $\xi$, and standard i.i.d.\ large deviation estimates
for the probabilities involving $z$.

\subsubsection{Upper tail bound}
It remains to prove
\be
P\left( R(t)> tr+t\e\right) \le A\exp(-Ct).
\label{tldput}
\ee
The argument is similar.
By statement (\ref{tXinf1}) and Lemma \ref{tindexlm} applied
to $x=r+\e$, we get
\beas
&&P\left( R(t)> tr+t\e\right) \\
&\le&
P\left( \mbox{$z_{[t(r+\e)]}(t)> z_k(0)-\xi^k_{[t(r-\e)]-k}(t)$
for all $k\ge 0$}
\right) \\
&\le&
P\left( \mbox{$z_{[t(r+\e)]}(t)= z_k(0)-\xi^k_{[t(r-\e)]-k}(t)$
for some $k< 0$}
\right) \\
&\le& e^{-Ct} +
P\left( \mbox{$z_{[t(r+\e)]}(t)= z_k(0)-\xi^k_{[t(r+\e)]-k}(t)$
for some $at\le k<0$}
\right).
\eeas
Check that $\rho y-g(r+\e-y)$ is strictly increasing
for $y\le\e$. Choose $\delta>0$ so that
\be
\mbox{$\rho\e-g(r) \ge \rho y-g(r+\e-y) +5\delta$ for all $y\le 0$.}
\label{tdeltadef2}
\ee
Choose a partition $a=y_0<y_1<\cdots<y_n=0$ of $[a,0]$ so that
$y_{i+1}-y_i<\delta$ for all $i$. Then
\be
\mbox{$t\rho\e-tg(r) \ge t\rho y_i-tg(r+\e-y_{i-1}) +4\delta t$
for all $t> 0$, $1\le i\le n$.}
\label{tineq3}
\ee
Reasoning as we did for the lower tail,
\beas
&&P\left( R(t)> tr+t\e\right) \\
&\le& e^{-Ct} +
\sum_{i=1}^n
P\left( z_{[t\e]}(0)-\xi^{[t\e]}_{[t(r+\e)]-[t\e]}(t)
\le z_{[ty_i]}(0)-\xi^{[ty_{i-1}]}_{[t(r+\e)]-[ty_{i-1}]}(t)
\right)\\
&\le& e^{-Ct} +
\sum_{i=1}^n \left\{
P\left( z_{[t\e]}(0)\le t\rho\e-\delta t\right)
+ P\left(
\xi^{[t\e]}_{[t(r+\e)]-[t\e]}(t)\ge tg(r)+\delta t\right)\right.\\
&&\qquad\left.
  + P\left(  z_{[ty_i]}(0)\ge  t\rho y_i+\delta t\right)
+ P\left(\xi^{[ty_{i-1}]}_{[t(r+\e)]-[ty_{i-1}]}(t) \le
tg(r+\e-y_{i-1}) -\delta t\right)
\right\}\\
&\le& Ae^{-Ct}.
\eeas
This completes the proof of Theorem \ref{secondclassestimate}.

\sect{Proofs of Theorems \ref{secondclass}, \ref{addfctvar}, and Corollary 
\ref{invarprin}}

The strategy of proof of the main theorems is to use the second-class
estimate and relation between occupation-times and second-class
particles, Theorem
\ref{secondclassestimate} and (\ref{relation}),  to establish
Theorem
\ref{addfctvar} for the particular totally asymmetric exclusion process in 
$d=1$
where $p(1)=1$ and $p(i)=0$ for $i\neq 1$.  Then, quoting a variance
comparison result (Proposition \ref{3.7} below), we generalize the
particular case to the full statement of Theorem \ref{addfctvar}.  Then, using
again the relation between occupation-times and
second-class
particles (\ref{relation}), we get Theorem
\ref{secondclass}.
Corollary \ref{invarprin} follows
as an easy consequence.

We will need a few preliminary results proved in \cite{FF} , \cite{S} and 
\cite{comparison}.
For $I=(i_1,\ldots,i_k) \subset \ZZ^d$ composed
of distinct vertices and
$k\geq 1$, define
``centered'' and ``monotone'' $k$-point functions respectively as
$$C^\rho_I (\eta) = (\eta_{i_1} - \rho)(\eta_{i_2}-\rho)
\cdots (\eta_{i_k} - \rho) \ \ {\rm and }$$
$$M^\rho_I(\eta) = (\eta_{i_1}\eta_{i_2}\cdots \eta_{i_k})- \rho^k.$$
It is useful to note that the monotone functions, $M^\rho_I(\eta)$ are
increasing local functions of $\eta$.
Observe now, for any local function $f(\eta)$, that
there exists $K=K_f<\infty$
such that $f$
can be represented in terms of a
finite linear combination of centered or monotone 
functions,
\begin{eqnarray*}
f &=& E_\rho[f] + \sum_{k=1}^K\sum_{|I|=k} \alpha_I C^\rho_I\\
&=& E_\rho[f] + \sum_{k=1}^K\sum_{|I|=k}\beta_I M^\rho_I
\end{eqnarray*}
with respect to some constants $\alpha_I$ and $\beta_I$.

\begin{lemma}
\label{3.1}
For any exclusion process with finite-range jump
rates $p$ and any $\rho\in [0,1]$ we have the following variance estimates.
There exists a constant $D_1=D_1(\rho, p)$ such that for
all $I\subset \ZZ^d$ such that $|I|=k$ we have
$$\sigma^2_t(C^\rho_I,\rho,p) \leq D_1t$$
when $k\geq 3$ in $d=1$, $k\geq 2$ in $d=2$, and $k\geq 1$ in $d\geq
3$.
Also, there exist constants $D_2=D_2(\rho,p)$ and $D_3=D_3(\rho,p)$
such that for all $i\in \ZZ^d$,
  $$\sigma^2_t(C^\rho_i - C^\rho_0)\leq D_2t,$$
and for all $i,j\in \ZZ^d$
$$\sigma^2_t(C^\rho_{(ij)} - C^\rho_{(01)})\leq D_3t.$$
\end{lemma}

{\it Proof.}
The proof follows directly from Lemma 3.9 \cite{S} (which bounds
$\sigma^2(f,\rho,p)_t \leq 10t\|f\|_{-1}(\rho,\bar{p})$ for local $f$)
and Lemma 3.4 \cite{S} (which bounds 
$$\|f\|_{-1}(\rho,\bar{p})<\infty \Leftrightarrow
\left\{ \begin{array}{ll}
E_\rho[f],\sum_{|I|=1}\alpha_I, \sum_{|I|=2}\alpha_I =0 &
{\rm when \ }d=1\\
E_\rho[f],\sum_{|I|=1}\alpha_I=0 & {\rm when\ } d=2\\
E_\rho[f]=0&{\rm when \ } d\geq 3,
\end{array}
\right.
$$
in terms of the centered basis representation).
\qed

Evidently, from this lemma, and the inequality 
$(a+b)^2 \leq 2a^2 + 2b^2$, the only
variances of centered functions
not bounded in $d=1$ are those of
$C^\rho_0$ and $C^\rho_{(01)}$, and in $d=2$ of
$C^\rho_0$.
The next lemma gives a relation between the two functions in $d=1$.

\begin{lemma}
\label{3.2}
We have in $d=1$ that
$$(\eta_0-\rho)(\eta_1-\rho) = [\rho(1-\rho)-\eta_0(1-\eta_1)]
+(1-2\rho)(\eta_0-\rho) + \rho[(\eta_0-\rho)-(\eta_1-\rho)].$$
\end{lemma}

{\it Proof.}
This follows from easy algebra. \qed

In $d=1$, the function $c(\eta)=\eta_0(1-\eta_1)$ arises in the study of 
the particle
current across the bond $0-1$.  Let $N_{0,1}(t)$ be the number of
particles which cross from $0$ to $1$ in time $t$.  Then $N_{0,1}(t)$
is a counting process with compensator $A_{c}(t) = \int_0^t
p(1)\eta_0(s)(1-\eta_1(s))ds$ so that
$M_{0,1}(t) = N_{0,1}(t) - A_{c}(t)$ is a square integrable martingale
with $E_\rho[M^2_{0,1}(t)] = p(1)\rho(1-\rho)t$.  The current has been
intensively studied in \cite{FF}.

\begin{lemma}
\label{3.3}
For the totally asymmetric nearest-neighbor exclusion process in $d=1$
with jump rate $p$, $p(1)=1$ and $p(i)=0$ for $i\neq 1$, we have that
$$\lim_{t\rightarrow \infty}
\frac{1}{t}E_\rho[(N_{0,1}(t) - p(1)\rho(1-\rho)t)^2] =
\rho(1-\rho)|1-2\rho|$$
and so, for all large $t$,
$$\sigma^2_t(c(\eta)-\rho(1-\rho), \rho, p)
\leq
3\rho(1-\rho)[1+|1-2\rho|]t.
$$
\end{lemma}

{\it Proof.}
The variance of $N_{0,1}(t)$ is
explicitly computed in
Theorem 1 \cite{FF}.
With the variance bound and the square
  martingale estimate, the
  inequality $(a+b)^2 \leq 2a^2 + 2b^2$ gives the 
last line.

  Alternatively, one can bypass the careful
computation
in \cite{FF} by observing that $N_{0,1}(t)$ has negatively correlated
increments and therefore has variance on the order $O(t)$.  Indeed,
let $N_{1,0}(t)$ be the number of particles crossing from $1$ to $0$
in time $t$, and write
\begin{eqnarray*}
&&E_\rho[(N_{0,1}(t)-\rho(1-\rho)t)(N_{0,1}(t+s) - N_{0,1}(t)-\rho(1-\rho)s)]\\
&&\ \ \ \ = E_\rho[(N_{0,1}(t)-\rho(1-\rho)t)E_{\eta(t)}[N_{0,1}(s)-\rho(1-\rho)s]]\\
&&\ \ \ \ = \int E^*_\eta[N_{1,0}(t)-\rho(1-\rho)t]E_\eta[N_{0,1}(s)-\rho(1-\rho)s] dP_\rho
\end{eqnarray*}
by time-reversal at time $t$ in the last line where $*$ refers to the
reversed process (for which also $P_\rho$ is invariant).  

Now, the
functions
$\phi(\eta) =  E^*_\eta[N_{1,0}(t)-\rho(1-\rho)t]$ and 
$\psi(\eta) = E_\eta[N_{0,1}(s)-\rho(1-\rho)s]$ have opposite
monotonicities.  That is, suppose 
$\eta$ and $\eta'$ are two configurations
such that $\eta_i = \eta_i'$ for all $i\neq x$ and $\eta_x =
1-\eta_x'=0$ for an $x\leq 0$.  By the basic coupling, the extra
particle at $x$ in the $\eta'$ configuration is a second-class
particle.  Let $N_{0,1}'(t)$ be the number of particles crossing from
$0$ to $1$ in time $t$ for the process begun at $\eta'$.
A moment's thought now
convinces that when the second-class particle is to
the left of $0$ or to the right of $1$ at time $t$, the
numbers $N_{0,1}(t)=N_{0,1}'(t)$ and $N_{0,1}(t) =N_{0,1}'(t)+1$
respectively.  Therefore, 
$\psi$ increases if $\eta$ is increased to the left of $0$.
Also, putting the extra particle initially
at $x\geq 1$ gives by an analogous argument that
$\psi$ decreases if $\eta$ is increased to the right of
$1$.  Similarly, 
as the jump rates are reversed in the adjoint
process, we have that $\phi$ decreases (increases) when 
$\eta$ is increased to the left of $0$ (increased to the right of $1$).

Finally, as $P_\rho$ is product measure, and therefore FKG, we have that $\int
\phi(\eta)\psi(\eta) dP_\rho \leq E_\rho[\phi]E_\rho[\psi]=0$. \qed
As a consequence, we have the following statement.
\begin{lemma}
\label{3.4}
For the totally asymmetric nearest-neighbor exclusion processes in $d=1$ 
with jump rate $p$, $p(1)=1$ and $p(i)=0$ for $i\neq 1$,
we have when
$\rho\neq 1/2$
that $\sigma^2_t(C^\rho_{01},\rho,p)\leq D_1t$
if $\sigma^2_t(C^\rho_0,\rho,p)\leq D_2t$ for some constants $D_1,
D_2$.  When $\rho=1/2$, already
$\sigma^2(C^\rho_{01},\rho,p)\leq D_3t$ for some constant
$D_3$.
\end{lemma}

{\it Proof.}
The bounds for $\rho\neq 1/2$ and $\rho=1/2$ follow
directly from Lemmas \ref{3.2} and \ref{3.3}.
\qed

We now state as a proposition the consequence of Theorem 
\ref{secondclassestimate} using the
relation (\ref{relation}).
\begin{proposition}
\label{3.5}
For the totally asymmetric nearest-neighbor
exclusion process in $d=1$ with jump rate $p(1)=1$ and $p(i)=0$ for
$i\neq 1$, we have when
$\rho\neq 1/2$ that
$$\sigma_t^2(C^\rho_0,\rho,p) \leq D t$$
for some constant $D=D(\rho)$.
\end{proposition}

{\it Proof.}
 From Theorem \ref{secondclassestimate} we have that the second-class
particle is $\P$-transient when $\rho\neq 1/2$ ($\Leftrightarrow v(\rho, p)\neq
0$) in
$d=1$.  Therefore, from (\ref{relation}), we have that
$$\sigma^2(C^\rho_0,\rho,p) = \lim_{t\rightarrow \infty}
t^{-1}\sigma^2_t(C^\rho_0,\rho,p)<\infty$$
exists and is finite. \qed

One of the results from \cite{S} is now quoted.
\begin{lemma}
\label{3.6}
For exclusion processes in $d\geq 1$ with finite-range
jump rates $p$, we have that $\sigma^2(f,\rho,p)$ exists whenever
$f$
is an increasing mean-zero $L^2(P_\rho)$ function.  In addition, if $f=f_+-f_-$ is the difference of two local
increasing mean-zero
functions whose limiting variances are finite, 
$\sigma^2(f_\pm,\rho,p)<\infty$, then also
$\sigma^2(f,\rho,p)<\infty$ exists.
\end{lemma}

{\it Proof.}
This follows from 
Lemma 3.1 \cite{S} (which gives
existence of $\sigma^2(f,\rho,p)$ when $f$ is (non-trivial) increasing, mean-zero, and in $L^2(P_\rho)$) 
and Lemma 3.2 \cite{S} (which proves
existence of $\sigma^2(f,\rho, p)<\infty$ when
$f=f_+-f_-$ for $f_+$ and $f_-$ which are local, 
increasing, mean-zero, and satisfy
$\sigma^2(f_\pm,\rho,p)<\infty$).
\qed

An application of the results in \cite{comparison} to limiting variances is the following.
\begin{proposition}
\label{3.7}
Consider the exclusion process in $d\geq 1$ with finite-range
jump rates $p$.  Define the nearest-neighbor
jump rate
$p'$, in terms of $p$, by
$$p'(\pm e^l)=
\left\{\begin{array}{ll}
                 \max[\pm\sum_i (i\cdot e^l)p(i),0]  & {\rm when \ }\sum_i
                 (i\cdot e^l)p(i)\neq 0 \\
                  1 & {\rm when \ }\sum_i (i\cdot e^l)p(i)=0
                 \end{array}
            \right. \\
$$
where $e^l$ for $1\leq l\leq d$ are the unit basis vectors in $\ZZ^d$
corresponding to the positive axes.
With respect to
the exclusion model corresponding to $p'$, we have,
when $f$ is a local increasing mean-zero function, that
$\sigma^2(f,\rho,p)<\infty$ if and only if
$\sigma^2(f,\rho,p')<\infty$.
\end{proposition}

{\it Proof.}
This is Corollary 6.1 of \cite{comparison}.
\qed
Note in $d=1$ that $p'$ reduces to a totally asymmetric
nearest-neighbor jump rate when
$p$ has non-zero drift, and to a symmetric one when $p$ is mean-zero.

The following lemma states one of the main
weak convergence results in $d=1,2$ found in \cite{S}.
\begin{proposition}
\label{3.8}
Consider
exclusion processes with finite-range jump
rates $p$ with non-zero drift in $d=1,2$.
Suppose that $f=f_+-f_-$ is the difference of two increasing local
mean-zero functions
such that $\sigma^2(f_\pm,\rho,p)<\infty$ so that, by Lemma \ref{3.6},
$\sigma^2(f,\rho,p)<\infty$ exists.
Then,
with respect to initial configurations given by
$P_\rho$, we have the weak convergence in $C[0,\infty)$
to Brownian motion $B$,
$$\lim_{\alpha\rightarrow \infty}
\frac{1}{\sqrt{\alpha}} A_f(\alpha t) =
B(\sigma^2(f,\rho,p) t).$$
\end{proposition}

{\it Proof.}
This is part (i) of Theorem 1.1 of \cite{S} for which 
the
invariance principle is proved
directly in the uniform topology when $p$ has
non-zero drift.  
\qed
\vskip .5cm
We now prove the main results.

{\bf Proof of Theorem {\ref{addfctvar}}.}
Let
  $f$ be a mean-zero local function, $E_\rho[f]=0$.
Then, by
the monotone basis expansion,
$f$ can be decomposed into the difference of two
local increasing mean-zero
functions, $f_+$ and $f_-$, where
$$\begin{array}{rll}
f_+ &=& \sum_{k=1}^K\sum_{\beta_I\geq 0} \beta_I M^\rho_I, \ \ {\rm and}\\
f_- &=&
\sum_{k=1}^K\sum_{\beta_I< 0} \beta_I M^\rho_I.
\end{array}
$$

We first show
the theorem for the $d=1$ 
totally asymmetric nearest-neighbor
model with $p(1)=1$ and $p(i)=0$ for $i\neq 1$.  In this case, when $\rho\neq 1/2$,
\be
\label{ordert}
\sigma^2_t(C^\rho_I,\rho,p) =
O(t)\ee
  for all sets $I$ from Lemmas \ref{3.1} and \ref{3.4},
and Proposition \ref{3.5}.  
 From Lemma \ref{3.6}, the limits $\sigma^2(f_\pm,\rho,p)$ both exist, and
from (\ref{ordert}) and repeated use of
the inequality $(a+b)^2\leq 2a^2 +2b^2$, they
  are both finite.  We can now
apply Lemma \ref{3.6} again to get
the statement of the theorem in this case.

Note that the theorem also holds in the totally asymmetric
nearest-neighbor model when the
jumps are to the left instead of right, or when the jump rate is
different from unity, by reflection, or time-change arguments respectively.

We now consider the general finite-range model with 
jump rate
$p$ with non-zero drift in $d=1$.  Observe that
$\sigma^2(h,\rho,p)$ exists for all local 
increasing mean-zero 
$h$ by Lemma \ref{3.6}, and
that
$$\sigma^2(h,\rho,p)<\infty  \ \Leftrightarrow \
\sigma^2(h,\rho,p')<\infty$$
by Proposition \ref{3.7}.  
As remarked after Proposition \ref{3.7},
$p'$ in $d=1$ is a totally asymmetric
nearest-neighbor jump rate.  For such rates, we have just
proved that
$\sigma^2(h,\rho,p')<\infty$ when $\rho\neq 1/2$.
In particular, 
we conclude $\sigma^2(f_\pm,\rho,p)<\infty$.
The full
theorem follows now, as before for the totally asymmetric
nearest-neighbor case, by invoking Lemma \ref{3.6}.
\qed
\vskip .5cm
{\bf Proof of Corollary \ref{invarprin}.}
This follows directly from Theorem \ref{addfctvar} and Proposition \ref{3.8}. \qed
\vskip .5cm
{\bf Proof of Theorem \ref{secondclass}.}
The $\P$-transience when $\rho\neq 1/2$ for $d=1$ exclusion models
with finite-range jump rates with non-zero drift follows from the
occupation-time relation (\ref{relation}) and the fact that
$\sigma^2(\eta_0-\rho,\rho,p)< \infty$ (Theorem \ref{addfctvar}). \qed

\bibliographystyle{plain}

\end{document}